# Experimental Results of the Search for Unitals in the Projective Planes of Order 25

## Stoicho D. Stoichev

Department of Computer Systems, Technical University of Sofia
email:stoi@tu-sofia.bg

**Abstract**
In this paper we present the results from a program developed by the author that finds the unitals of the known 193 projective planes of order 25. There are several planes for which we have not found any unital. One or more than one unitals have been found for most of the planes. The found unitals for a given plane are nonisomorphic each other. There are a few unitals isomorphic to a unital of another plane.
A t - (v; k; λ) *design* D is a set X of points together with a family B of k-subsets of X called blocks with the property that every t points are contained in exactly λ blocks. The design with t = 2 is called a *block-design.* The block-design is symmetric if the role of the points and blocks can be changed and the resulting confguration is still a block-design. A projective plane of order n is a symmetric 2-design with $v = n^2 + n + 1$, $k = n + 1$, $λ = 1$. The blocks of such a design are called lines. A unital in a projective plane of order $n = q^2$ is a set U of $q^3 + 1$ points that meet every line in one or q + 1 points. In the case projective planes of order n = 25 we have: q = 5 , the projective plane is 2 - (651; 26; 1) design, the unital is a subset of $q^3 + 1 = 5^3 + 1 = 126$ points and every line meets 1 or 6 points from the subset

## Article Outline

## 1. Introduction

We assume familiarity with the basics of the combinatorial designs (cf., e.g. [1]). A t - (v; k; λ) *design* D [1] is a set X of points together with a family B of k-subsets of X called blocks with the property that every t points are contained in exactly λ blocks. The design with t = 2 is called a *block-design.* The block-design is symmetric if the role of the points and blocks can be changed and the resulting confguration is still a block-design. A projective plane of order n is a symmetric 2-design with $v = n^2 + n + 1$, $k = n + 1$, $λ = 1$. The blocks of such a design are called lines. A unital in a projective plane of order $n = q^2$ is a set U of $q^3 + 1$ points that meet every line in one or q + 1 points.
In the case of projective planes of order n = 16 we have: q = 4 , the projective plane is 2 - (273; 17; 1) design, the unital is a subset of $q^3 + 1 = 4^3 + 1 = 65$ points and every line meets 1 or 5 points from the subset.
In the case of projective planes of order n = 25 we have: q = 5 , the projective plane is 2 - (651; 26; 1) design, the unital is a subset of $q^3 + 1 = 5^3 + 1 = 126$ points and every line meets 1 or 6 points from the subset.

The experimental results of the algorithm in [5] for the known planes of order 16 [4] are given in [6].

## 2. Experimental results

In [8] we present the text of the current paper and the experimental results (in Appendix) from a program that finds the unitals in all 193 known projective planes of order 25 from the website [3].. The vertex labels in this websites start from 0, but in our program and results the starting label is 1..The program we use in present paper is based on the algorithm described in [5]:. In [5] we describe the algorithm for finding unitals and maximal arcs in projective planes of order 16. In this paper we use the same algorithm but with parameters for a plane of order 25.  The algorithm is heuristic – it does not guarantee finding all possible unitals of a given plane. All fond unitals of an unital are nonisomorphic each other. There are planes with no found unitals.

The results for each plane contain values for the following variables:
1. Name of the plane
2. Order of the plane automorphism group
3. Number of the orbits of the plane automorphism gtoup
4. Order of the unital automorphism group
5. Number of the orbits of the unital automorphism group
6. Sizes of the orbits of the unital automorphism group
7. Labels of the vertices of the unital

The information for 1 to 3 is present for each plane and the information for 4 to 7 is present for each found unital.

***Example*** (for the plane A1.HTM, see the text below - from the line with 'A1.HTM' : The order of the plane automorphism group is 360000 and the number of its orbits is 8. Then, the results for all 9 found unitals follow. The orders of their automorphism groups are 144, 24, 20, 20, 15, 10, 10, 10, 10.

```
A1.HTM
ORDER OF THE PLANE AUTOMORPHISM GROUP              360000
NUMBER OF THE ORBITS OF THE PLANE AUTOMORPHISM GTOUP=    8
ORDER OF THE UNITAL AUTOMORPHISM GROUP=                 144
NUMBER OF THE ORBITS OF THE UNITAL AUTOMORPHISM GROUP =           3
SIZES OF THE ORBITS OF THE UNITAL AUTOMORPHISM GROUP=   1-  72    2-  48
 3-   6
UNITAL=
      6    7    8    9   10   11   30   31   32   34   35   39  181  182  183
 184
    185  186  187  188  189  190  191  192  193  194  196  197  198  199  208
 209
    210  211  212  213  214  215  216  229  232  237  238  239  240  241  242
 243
    244  245  246  247  248  249  250  251  252  266  270  271  274  275  276
 277
    278  279  289  290  291  292  293  294  295  296  297  298  299  300  301
 303
    304  307  308  309  353  355  357  359  361  363  364  365  366  367  368
 369
    376  378  379  381  384  387  606  607  608  609  611  615  616  617  618
 619
    620  621  628  632  634  635  636  638  640  643  646  647  648  651
ORDER OF THE UNITAL AUTOMORPHISM GROUP=                  24
NUMBER OF THE ORBITS OF THE UNITAL AUTOMORPHISM GROUP =           8
SIZES OF THE ORBITS OF THE UNITAL AUTOMORPHISM GROUP=   1-  24    2-  24
 3-  24    4-  24    5-  24    6-   4    7-   1    8-   1
UNITAL=
```

```
     1    10    11    22    23    25    40    41    42    49    50    51    64    65    67    68
    70    71    85    91    94    95    97    98   118   119   120   125   127   133   134   138
   142   148   149   154   173   174   176   177   178   179   183   184   186   190   195   196
   197   200   202   205   206   207   210   214   216   223   225   226   227   229   231   233
   238   239   245   249   251   254   259   262   264   267   271   273   275   278   326   329
   331   336   340   345   389   392   397   398   406   409   410   411   412   413   420   423
   431   435   442   446   451   453   460   462   480   481   489   492   538   541   543   553
   570   573   576   577   593   595   600   603   608   611   619   620   635   636
ORDER OF THE UNITAL AUTOMORPHISM GROUP=                   20
NUMBER OF THE ORBITS OF THE UNITAL AUTOMORPHISM GROUP =           8
SIZES OF THE ORBITS OF THE UNITAL AUTOMORPHISM GROUP=   1- 20   2- 20
  3- 20   4- 20   5- 20   6- 20   7-  5   8-  1
UNITAL=
     8    27    29    30    32    48    49    53    61    62    66    71    79    82   101   104
   110   115   121   129   133   138   143   149   159   162   168   173   195   197   202   204
   206   207   220   225   229   237   239   241   243   245   248   257   263   264   266   276
   281   286   289   293   297   306   315   319   325   327   329   332   334   337   339   345
   354   358   367   370   377   385   392   396   405   409   411   415   419   423   424   431
   437   446   450   453   456   458   470   472   473   477   478   489   490   495   498   505
   515   518   519   523   524   531   533   536   537   544   554   557   559   560   567   568
   571   573   587   595   602   603   604   605   609   619   620   626   640   641
ORDER OF THE UNITAL AUTOMORPHISM GROUP=                   20
NUMBER OF THE ORBITS OF THE UNITAL AUTOMORPHISM GROUP =           8
SIZES OF THE ORBITS OF THE UNITAL AUTOMORPHISM GROUP=   1- 20   2- 20
  3- 20   4- 20   5- 20   6- 20   7-  5   8-  1
UNITAL=
     8    27    32    33    35    41    51    53    55    61    71    72    74    89    91    97
    98   101   104   111   117   119   121   133   140   158   167   168   179   195   205   206
   208   216   217   220   225   230   236   241   245   254   256   259   260   266   270   276
   280   282   290   293   295   297   302   314   315   329   330   339   340   341   347   350
   360   362   366   368   371   375   377   378   401   405   410   411   414   421   425   439
   440   446   449   456   458   460   476   478   482   486   488   492   498   512   514   523
   527   531   533   534   536   537   540   542   544   549   551   554   555   557   560   567
   568   573   574   603   604   613   621   624   629   630   634   637   643   648
ORDER OF THE UNITAL AUTOMORPHISM GROUP=                   15
NUMBER OF THE ORBITS OF THE UNITAL AUTOMORPHISM GROUP =          10
SIZES OF THE ORBITS OF THE UNITAL AUTOMORPHISM GROUP=   1- 15   2- 15
  3- 15   4- 15   5- 15   6- 15   7- 15   8- 15   9-  5  10-  1
UNITAL=
```

```
    24    28    36    45    48    49    51    56    57    66    67    68    75    79    86
    88
    98   107   117   124   125   127   134   135   139   140   143   146   150   152   157
   160
   161   163   164   173   177   179   195   200   202   204   206   212   215   232   242
   244
   245   270   275   280   283   288   293   295   296   297   303   309   315   318   319
   323
   326   327   332   333   348   361   364   365   371   372   373   379   380   384   385
   389
   391   403   406   411   416   420   427   435   437   438   441   444   448   451   454
   464
   466   470   475   482   485   488   491   495   537   555   556   562   564   581   598
   599
   605   609   615   618   622   623   626   627   635   636   643   644   647   650
ORDER OF THE UNITAL AUTOMORPHISM GROUP=                    10
NUMBER OF THE ORBITS OF THE UNITAL AUTOMORPHISM GROUP =              14
SIZES OF THE ORBITS OF THE UNITAL AUTOMORPHISM GROUP=   1-  10   2-  10
  3-  10   4-  10   5-  10   6-  10   7-  10   8-  10   9-  10  10-  10  11-
 10  12-  10  13-   5  14-   1
UNITAL=
    10    27    31    37    39    43    46    54    56    61    64    67    72    78    85
    88
    95   100   114   116   118   131   133   147   148   153   155   163   174   182   190
   201
   209   212   214   216   217   218   226   234   244   246   254   256   258   260   263
   264
   283   285   287   294   302   303   306   310   315   320   326   327   343   346   349
   350
   356   359   363   364   370   373   374   376   394   398   401   411   413   418   424
   425
   433   443   448   452   458   463   465   467   470   475   479   495   504   511   516
   518
   524   530   536   539   541   544   546   549   553   555   564   565   568   570   572
   593
   595   597   601   603   606   608   610   622   629   630   634   639   645   648
ORDER OF THE UNITAL AUTOMORPHISM GROUP=                    10
NUMBER OF THE ORBITS OF THE UNITAL AUTOMORPHISM GROUP =              14
SIZES OF THE ORBITS OF THE UNITAL AUTOMORPHISM GROUP=   1-  10   2-  10
  3-  10   4-  10   5-  10   6-  10   7-  10   8-  10   9-  10  10-  10  11-
 10  12-  10  13-   5  14-   1
UNITAL=
    10    27    30    35    36    41    48    57    58    59    69    75    78    81    83
    85
    87    91    99   103   106   114   117   131   138   145   148   153   158   161   162
   171
   176   184   188   191   192   197   200   211   216   220   222   223   230   233   251
   254
   259   268   269   276   278   283   284   291   304   307   309   311   319   322   325
   341
   359   361   366   379   380   384   388   396   398   403   409   421   422   423   429
   441
   442   447   449   450   451   457   469   473   474   476   485   488   491   493   499
   509
   510   512   513   520   528   530   535   539   544   550   551   564   568   574   576
   580
   581   584   594   595   597   603   609   613   614   615   620   627   642   643
ORDER OF THE UNITAL AUTOMORPHISM GROUP=                    10
NUMBER OF THE ORBITS OF THE UNITAL AUTOMORPHISM GROUP =              14
```

```
SIZES OF THE ORBITS OF THE UNITAL AUTOMORPHISM GROUP=   1- 10    2- 10
 3- 10    4- 10    5- 10    6- 10    7- 10    8- 10    9- 10   10- 10   11-
10   12- 10   13-  5   14-  1
UNITAL=
    10    27    55    56    60    64    74    78    82    85    89    94    95    98   100
108
   114   117   120   126   131   132   133   135   143   148   151   153   161   174   175
179
   180   193   201   205   206   207   208   215   216   223   226   234   235   238   243
245
   249   254   255   259   261   263   266   267   283   285   292   294   304   307   308
310
   314   318   320   328   331   334   335   340   344   346   356   358   359   370   377
378
   388   394   395   397   398   400   432   435   437   439   455   458   467   479   481
482
   491   493   501   505   506   509   521   530   539   544   548   556   561   564   568
585
   588   595   596   597   598   603   606   610   622   632   638   639   647   649
ORDER OF THE UNITAL AUTOMORPHISM GROUP=                      10
NUMBER OF THE ORBITS OF THE UNITAL AUTOMORPHISM GROUP =              14
SIZES OF THE ORBITS OF THE UNITAL AUTOMORPHISM GROUP=   1- 10    2- 10
 3- 10    4- 10    5- 10    6- 10    7- 10    8- 10    9- 10   10- 10   11-
10   12- 10   13-  5   14-  1
UNITAL=
    10    27    28    30    35    41    45    56    62    64    71    78    84    85    87
 95
   100   106   107   114   119   124   131   133   138   148   153   166   167   174   177
181
   186   189   191   194   201   216   220   225   226   231   234   251   254   263   268
271
   274   280   281   283   285   289   294   301   310   317   319   320   324   341   342
346
   354   355   356   359   361   366   370   380   389   394   398   403   405   409   423
445
   450   457   458   460   466   467   472   476   479   484   485   492   499   500   514
515
   519   530   531   535   537   538   539   544   550   564   568   576   582   584   595
597
   600   602   603   606   610   613   616   622   626   627   639   641   642   650
```

## 3. Concluding remarks

By our algorithm we have found new unitals in projective plane of order 25, but not all of them.
The following approaches can be used to find more or all unitals : (a) Development of improved algorithms by finding new conditions for pruning the search tree;
(b) Transformation of the solution for one plane to solution for another plane (R. Mathon's approach - in private communication); (c) Development of parallel algorithms.

### **Acknowledgements**

The author would like to thank Vladimir Tonchev for suggesting the problem of developing an algorithm for finding unitals in projective planes and for giving the general idea of such an algorithm - use of unions of orbits, and for extensive discussions and exchanges for many years.
.